\documentclass[letterpaper, 10 pt, conference]{ieeeconf}  

\IEEEoverridecommandlockouts                              
\usepackage{amsmath, amssymb, amsxtra, mathrsfs}
\usepackage{graphics, epsfig, here, epstopdf, fancyhdr,graphicx}
\usepackage{bbm}
\usepackage{color}
\usepackage{multirow}
\usepackage[all]{xy}

\newcommand{\tS}{\textstyle}

\newcommand{\pf}{\noindent \textit{Proof:}\quad}
\newcommand{\boX}{\hfill $\Box$}
\newcommand{\bN}{\mathbb{N}}

\newcommand{\bR}{\mathbb{R}}

\newcommand{\cP}{\mathcal{P}}
\newcommand{\cC}{\mathcal{C}}
\newcommand{\cE}{\mathcal{E}}

\newcommand{\cL}{\mathcal{L}}

\newcommand{\cU}{\mathcal{U}}
\newcommand{\cV}{\mathcal{V}}

\newcommand{\cX}{\mathcal{X}}

\newcommand{\fA}{\,\forall\,}

\newcommand{\dS}{\displaystyle}

\newcommand{\bq}[1]{{\left[#1\right]}}
\newcommand{\bp}[1]{{\left(#1\right)}}

\newcommand{\vphi}{\varphi}

\newcommand{\cG}{\mathcal{G}}

\newcommand{\eP}{\varepsilon}

\newcommand{\pp}{\partial}

\newcommand{\rd}{\mathrm{d}}
\newcommand{\ddx}[1]{{\frac{\rd }{\rd #1}}}
\newcommand{\ppx}[1]{{\frac{\partial }{\partial #1}}}

\newcommand{\pxpy}[2]{{\frac{\partial #1}{\partial #2}}}

\usepackage{algorithmic}
\usepackage{epstopdf}

\title{\LARGE \bf
Monotonicity of Actuated Flows on Dissipative Transport Networks
}

\author{Anatoly Zlotnik, Sidhant Misra, Marc Vuffray and Michael Chertkov
\thanks{A. Zlotnik, S. Misra, M. Vuffray, and M. Chertkov are with Los Alamos National Laboratory, Los Alamos, NM 87544.  \qquad
Email: \{azlotnik $\mid$ sidhant $\mid$ vuffray $\mid$ chertkov\}@lanl.gov.}}

\newtheorem{thm}{Theorem}

\newtheorem{prop}[thm]{Proposition}
\newtheorem{cor}[thm]{Corollary}

\newtheorem{df}{Definition}

\newtheorem{rem}{Remark}

\begin{document}

\maketitle
\thispagestyle{empty}
\pagestyle{empty}

\begin{abstract}
We derive a monotonicity property for general, transient flows of a commodity transferred throughout a network, where the flow is characterized by density and mass flux dynamics on the edges with density continuity and mass balance conditions at the nodes.  The dynamics on each edge are represented by a general system of partial differential equations that approximates subsonic compressible fluid flow with energy dissipation.  The transferred commodity may be injected or withdrawn at any of the nodes, and is propelled throughout the network by nodally located compressors.  These compressors are controllable actuators that provide a means to manipulate flows through the network, which we therefore consider as a control system.  A canonical problem requires compressor control protocols to be chosen such that time-varying nodal commodity withdrawal profiles are delivered and the density remains within strict limits while an economic or operational cost objective is optimized.  In this manuscript, we consider the situation where each nodal commodity withdrawal profile is uncertain, but is bounded within known maximum and minimum time-dependent limits.  We introduce the monotone parameterized control system property, and prove that general dynamic dissipative network flows possess this characteristic under certain conditions.  This property facilitates very efficient formulation of optimal control problems for such systems in which the solutions must be robust with respect to commodity withdrawal uncertainty. We discuss several applications in which such control problems arise and where monotonicity enables simplified characterization of system behavior.
\end{abstract}

\section{Introduction} \label{secintro}

The optimal allocation of commodity flows over networks has been studied from theoretical and computational perspectives since the early work of Ford and Fulkerson, which focused on maximal utilization of capacity and minimization of economic cost in the steady state \cite{ford58b,ford62}.  
Subsequent algorithms are prominent in operations research, with particular importance for transportation problems \cite{gass90}, which may involve commodities such as vehicles \cite{ma04,tjandra03phd,como10}, fluids \cite{wong68,misra15}, 
energy \cite{geidl07}, and information \cite{intanagonwiwat02,dapice08}.

The difficulty of network flow problems is amplified when the flows are unbalanced, i.e., when commodity inflows at origins and outflows at destinations are time-dependent \cite{gottlich05}.  Such situations arise in air traffic flow \cite{ma04}, telecommunications networks \cite{dapice08}, and other flow problems that require dynamic modeling and control design \cite{herty03,como10}.  The number of constraints and decision variables increases by a factor directly related to the temporal complexity of commodity inflows and outflows.  The dynamics are then characterized by systems of ordinary or partial differential equations (ODEs or PDEs), which represent fluid flow or the aggregated motion of discrete particles.  In this context, optimization requires incorporating differential constraints rather than purely algebraic ones, for example in vehicle traffic \cite{gugat05} and natural gas pipeline flows \cite{steinbach07pde}.

A further challenge to computational tractability for applications arises through the presence of uncertainty in the volume and timing of the variable commodity inflows and outflows.  It is desirable for control or routing policies to be feasible for any instance under such uncertainty \cite{como10}, and the notion of robust optimization has been applied to create solutions in the discrete case using global information \cite{bental98,bertsimas03}.  In the case of continuous dynamic flows, uncertainty in constant or time-varying functional system parameters, e.g., network inflows and outflows, requires a continuum of constraints to ensure feasibility of the optimization solution.  The challenge then becomes to show that feasibility for a finite number of appropriate scenarios will guarantee feasibility for an entire such uncountable ensemble of constraints.

A recent approach to control uncertain network flows sidesteps the need for global optimization over a possibly non-convex landscape by examining stability and robustness of distributed routing solutions \cite{como13a,como13b}.  The methodology in these studies was enabled by demonstrating that the dynamics in question were monotone control systems \cite{angeli03,como10,lovisari14}.  Such so-called cooperative systems, which possess monotonicity with respect to certain input variables, were investigated in the context of ordinary differential equation theory \cite{kamke32,hirsch85,smith88,hirsch05}.  Recently, monotonicity properties were found to facilitate analysis of chemical reaction networks, power systems, and turbulent jet flows \cite{deleenheer04,budivsic12}, and propagation of order properties for stochastic systems have been proposed \cite{sootla15}.  Monotonicity has also recently been established for steady-state commodity flows on networks in order to enable efficient algorithms for robust optimization of natural gas systems under uncertainty \cite{vuffray15cdc}. 

In this manuscript, we derive monotonicity properties for a class of actuated dynamic flows described by PDEs coupled at the boundaries in a network structure.  Under certain conditions, we show that commodity density anywhere in the network can only increase monotonically when any commodity injection is increased.  We first produce an ODE system by spatial discretization of the flow network using lumped elements.  The resulting model is described as a parameterized control system, to which we apply the standard Kamke conditions \cite{kamke32,hirsch05} in order to establish monotonicity with respect to parameter functions.  We also present conditions on local feedback policies that maintain monotonicity, and describe how robust optimal control problems for monotone systems can be compactly formulated.

The manuscript is organized as follows.  In Section \ref{sec:formulation}, we formulate a class of actuated commodity flows through dissipative transport networks as a system of PDEs over a collection of domains that form a graph when coupled by Kirchhoff-Neumann boundary conditions.  Section \ref{sec:discretization} describes a lumped-element spatial discretization of the continuum dynamics, in which the network is refined and nodal density dynamics are obtained.  Section \ref{sec:monotonicity}
establishes monotonicity of the nodal dynamics by applying well-known monotone systems results.  In Section \ref{sec:application} we discuss the applications to monotonicity-preserving local feedback control and robust optimal control, and  conclude with Section \ref{sec:conc}.

\vspace{-1ex}
\section{Dynamic Dissipative Flows on Networks} \label{sec:formulation}

We consider a network with flows of a compressible fluid commodity through pipelines that are connected at junctions where the fluid can be compressed into a pipe, or withdrawn from or injected into the system.  This network is represented as an oriented weighted graph $\Gamma=\left(\cV,\cE,\lambda\right)$ where $\cV$ is the set of vertices and $\cE\subset V\times V$ is the set of directed  edges $(i,j)\in\cE$ that connect the nodes $i,j\in\cV$.  The incoming neighborhood of $i\in \cV$ is denoted
by $\partial_{+}i=\left\{ j\in \cV\mid(j,i)\in \cE\right\} $
and its outgoing neighborhood is denoted by $\partial_{-}i=\left\{ j\in \cV\mid(i,j)\in \cE\right\}$.  Every edge $(i,j)\in \cE$ is associated with a spatial dimension on the interval $I_{ij}=[0,L_{ij}]$, where $L_{ij}=\lambda(i,j)$ and $\lambda:\cE\to\bR_+$ (where $\bR_+=\{x\in\bR\,:\,x\geq 0\}$) defines the graph edge weights corresponding to pipe lengths.  We let $V=|\cV|$ and $E=|\cE|$ denote the number of nodes and of edges, respectively.

The instantaneous state within each edge $(i,j)\in\cE$ is characterized by space-time dependent mass flux $\phi_{ij}:[0,T]\times I_{ij}\rightarrow\mathbb{R}$
and density $\rho_{ij}:[0,T]\times I_{ij}\rightarrow\mathbb{R}_{+}$ functions.  By convention, $\phi_{ij}(t,x_{ij})=-\phi_{ji}(t,L_{ij}-x_{ij})$.  In addition, every node $i\in \cV$ is associated with a time-dependent internal nodal density $\rho_{i}(t):[0,T]\rightarrow\mathbb{R}_{+}$ and is subject to
a time-dependent mass flux injection $q_{i}:[0,T]\rightarrow\mathbb{R}$.  We define a convention where $q_{i}$ is positive when the commodity is injected into the network at node $i\in\cV$, and is negative when the commodity is withdrawn.  

We suppose that the density and mass flux dynamics on the edge  $(i,j)\in\cE$ evolve according to the generalized dissipative distributed relations
\begin{align}
\dS \pp_t\rho_{ij}(t,x_{ij})+\pp_x\phi_{ij}(t,x_{ij}) & =  0 \label{eq:in_continuity} \\
\phi_{ij}(t,x_{ij})+f_{ij}(t,\rho_{ij}(t,x_{ij}), \partial_{x}\rho_{ij}(t,x_{ij})) & =0,  \label{eq:in_dissipation_eq}
\end{align}
which are called the continuity and dissipation equations.  

We define a set of controllable nodal actuators $\cC\subset \cE\times\{+,-\}$, where $(i,j)\equiv\{(i,j),+\}\in\cC$ denotes a controller located at node $i\in\cV$ that augments the density of the commodity that flows into edge $(i,j)\in\cE$ in the positive direction, while $(j,i)\equiv\{(i,j),-\}\in\cC$ denotes a controller located at node $j\in\cV$ that augments density into edge $(i,j)\in\cE$ in the negative direction.  Compression is then modeled as a multiplicative ratio $\underline{\alpha}_{ij}:[0,T]\to\bR_+$ for $\{(i,j),+\}\in\cC$ and $\overline{\alpha}_{ij}:[0,T]\to\bR_+$ for $\{(i,j),-\}\in\cC$.

\begin{figure}[t]
\centering
\includegraphics[width=.95\linewidth]{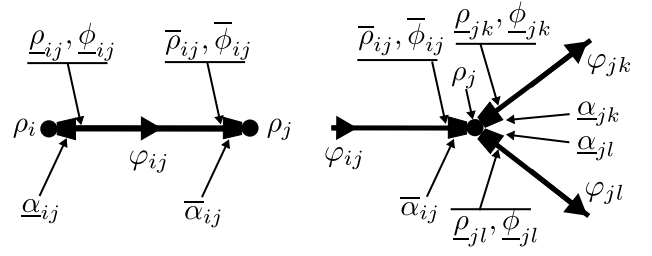} \vspace{-2ex} \caption{Nodal densities $\rho_j$, and edge endpoint variables $\underline{\rho}_{ij}$, $\underline{\phi}_{ij}$, $\overline{\rho}_{ij}$, and $\overline{\phi}_{ij}$, and actuators $\underline{\alpha}_{ij}$, $\overline{\alpha}_{ij}$ for an edge (left) and for a joint (right).} \label{fig:compatibility}
\end{figure}

Next, we establish nodal relations that characterize the boundary conditions for the flow dynamics \eqref{eq:in_continuity}-\eqref{eq:in_dissipation_eq} on each edge of the network.  For this purpose, we define
\begin{align}
\underline{\rho}_{ij}(t)\triangleq\rho_{ij}(t,0), \quad \overline{\rho}_{ij}(t)\triangleq\rho_{ij}(t,L_{ij}), \label{eq:end_p_def} \\
\underline{\phi}_{ij}(t)\triangleq\phi_{ij}(t,0), \quad \overline{\phi}_{ij}(t)\triangleq\phi_{ij}(t,L_{ij}), \label{eq:end_q_def}
\end{align}
and $\vphi_{ij}(t)\triangleq\phi_{ij}(t,\tS\frac{1}{2}L_{ij})$ in order to simplify notation.  At each vertex $i\in V$ the mass flux and density values at the endpoints of adjoining edges must satisfy certain compatibility conditions.  First, a Kirchhoff-Neumann property of mass conservation is ensured through nodal continuity equations
\begin{align}
q_j(t)+\sum_{i\in\partial_{+}j}\overline{\phi}_{ij}- \sum_{k\in\partial_{-}j}\underline{\phi}_{jk}=0, \quad \fA j\in\cV. \label{eq:in_nodal_continuity}
\end{align}
In addition, compatibility conditions for densities may involve jump discontinuities in space due to compression actuators, and are of the form
\begin{align}
\underline{\rho}_{ij}(t) = \underline{\alpha}_{ij}\rho_{i}(t), \quad
\overline{\rho}_{ij}(t)  = \overline{\alpha}_{ij}\rho_{j}(t), \quad \fA (i,j)\in\cE, \label{eq:in_pressure_comp}
\end{align}
where $\underline{\alpha}_{ij}$ and $\overline{\alpha}_{ij}$ are
positive compression ratios that represent actuation at the initial ($i\in\cV$) and endpoint ($j\in\cV$) vertex of each edge $(i,j)\in\cE$.  The variables $\rho_{i}$ are auxiliary variables that denote internal nodal density values.  The above compatibility conditions are visualized in Figure \ref{fig:compatibility}.  The instantaneous state of the system at time $t=0$ can be specified by initial density profiles of the form
\begin{align}
\!\!\! \rho_{ij}(0,x)=\rho_{ij}^{0}(x), \,\, \phi_{ij}(0,x)=\phi_{ij}^{0}(x), \quad \fA (i,j)\in\cE. \label{eq:in_initial_condition}
\end{align}

\begin{rem} \label{rem:wellposed} We suppose for some $k\in\bN$ that $\rho_{ij}^{0},\phi_{ij}^{0}\in C^k([0,L_{ij}])$, that \eqref{eq:in_nodal_continuity} and \eqref{eq:in_pressure_comp} hold at $t=0$, that the actuator functions satisfy $\underline{\alpha}_{ij},\overline{\alpha}_{ij}\in C_+^k([0,T])$ for $(i,j)\in\cE$, and that injection functions satisfy $q_i\in C^k([0,T])$ for $i\in\cV$.  It is assumed that the dynamics characterized by \eqref{eq:in_continuity}-\eqref{eq:in_dissipation_eq} with the initial conditions \eqref{eq:in_initial_condition} admit a unique classical solution on the interval $[0,T]$ under the above conditions.
\end{rem}


\section{Nodal Density Dynamics on Refined Network} \label{sec:discretization}

We use a lumped element approximation to characterize edge dynamics \eqref{eq:in_continuity} and \eqref{eq:in_dissipation_eq}, with nodal conditions \eqref{eq:in_nodal_continuity} and \eqref{eq:in_pressure_comp} and subject to injection profiles $q_{i}(t)$,  which approximately defines the state on the network in terms of nodal densities $\rho_j(t)$. Our approach is to add enough nodes to the network so that density and flow are nearly uniform on any given segment. 
In particular, we obtain dynamic equations where the state is represented by the vector of nodal densities $\rho=(\rho_1,\ldots,\rho_V)$.  We begin with the following definition.

\begin{df}[Spatial Graph Refinement] \label{def:graphref} The refinement $\hat{\cG}_{\eP}=(\hat{\cV}_{\eP},\hat{\cE}_{\eP},\hat{\lambda}_{\eP})$ of a weighted oriented graph $\cG=(\cV,\cE,\lambda)$ is made by adding nodes to $\cV$ to sub-divide edges of $\cE$ where the length $\hat{L}_{ij}\in\hat{\cL}_{\eP}$ of a new edge $(i,j)\in\hat{\cE}_{\eP}$ is satisfies
\begin{align}
\frac{\eP L_{\mu(ij)}}{\eP+L_{\mu(ij)}}<\hat{L}_{ij}<\eP,
\end{align}
where $\mu:\hat{\cE}\to\cE$ is an surjective map of refined edges to the parent edges in $\cE$.
\end{df}

\begin{rem} Spatial graph refinement preserves the network structure represented by the graph, and can finely discretize the coupled one-dimensional domains on which the network dynamics \eqref{eq:in_continuity}-\eqref{eq:in_dissipation_eq} with \eqref{eq:in_nodal_continuity}-\eqref{eq:in_pressure_comp} evolve.  For $\eP\ll\min_{{i,j}\in\cE} L_{ij}$, the lengths in $\hat{\cL}_{\eP}$ are nearly uniform and very close to $\eP$.
\end{rem}

\begin{rem} We assume that $\eP$ is small enough that the relative difference of density and flux at the start and end of each new edge $(i,j)\in\hat{\cE}_{\eP}$ is small.  Specifically,
\begin{align} \label{eq:pres_rel}
2\frac{\overline{\rho}_{ij}(t)-\underline{\rho}_{ij}(t)}{ \overline{\rho}_{ij}(t)+\underline{\rho}_{ij}(t)} \ll 1, \,\, 2\frac{\overline{\phi}_{ij}(t)-\underline{\phi}_{ij}(t)}{ \overline{\phi}_{ij}(t)+\underline{\phi}_{ij}(t)} \ll 1, \,\, \fA t
\end{align}
for the transient regime of interest.  In other words, $\eP$ is sufficiently small so that the relative density difference between neighboring nodes is very small at all times.
\end{rem}

\begin{figure}[t]
\centering{
\includegraphics[width=.95\linewidth]{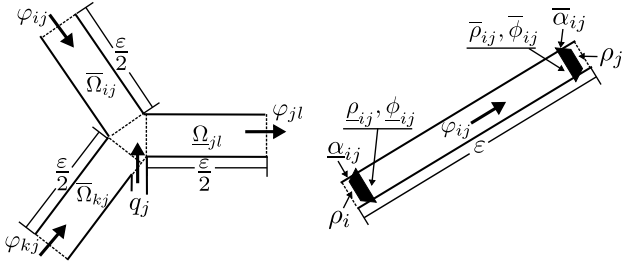} \caption{Lumped elements for discretization of equations \eqref{eq:in_continuity} at a node (left) and \eqref{eq:in_dissipation_eq} over an edge (right).} \label{fig:netconstit} }
\end{figure}

Consider an actuated flow network with a spatial graph refinement $\hat{\cG}_{\eP}$ with  $V_{\eP}=|\cV_{\eP}|$ nodes and $E_{\eP}=|\cE_{\eP}|$ edges of approximate length $\eP$.   Figure \ref{fig:netconstit} illustrates an example joint (left) and pipe (right).  The variable $q_j$ denotes an injection into the network at node $j$, and $\overline{\Omega}_{ij}$ and $\underline{\Omega}_{jk}$  are  sub-elements corresponding to halves of incoming and outgoing edges $(i,j)$ and $(j,k)$ in $\hat{\cE}_{\eP}$.  The flow at the midpoint of an edge is denoted $\vphi_{ij}=\phi_{ij}(t,\hat{L}_{ij}/2)$.   The densities $\underline{\rho}_{ij}$ and $\overline{\rho}_{ij}$ at the ends of the edge $(i,j)\in\cE$ are related to the nodal densities $\rho_i$ and $\rho_j$ by Equation \eqref{eq:in_pressure_comp}, as described in Section \ref{sec:formulation}.

First, we approximate the rate of change of mass within the nodal element by summing the integrals of mass flux gradient on each adjoining pipe segment.  That is,
\begin{align}
&\sum_{i\in\partial_+j}\int_{\overline{\Omega}_{ij}}\partial_x\phi_{ij}(t,x_{ij}) +\sum_{k\in\partial_-j}\int_{\underline{\Omega}_{jk}}\partial_x\phi_{jk}(t,x_{jk}) \label{eq:fx_constit1}  \\
&\qquad=  \sum_{i\in\partial_+j}(\overline{\phi}_{ij}-\vphi_{ij}) +\sum_{k\in\partial_-j}(\vphi_{jk}-\underline{\phi}_{jk})  \label{eq:fx_constit2}\\
&\qquad= \sum_{k\in\partial_-j}\vphi_{jk} - \sum_{i\in\partial_+j}\vphi_{ij} -q_j,
\label{eq:fx_constit3}
\end{align}
where the last step is due to the nodal balance condition \eqref{eq:in_nodal_continuity}. Next, applying mass conservation \eqref{eq:in_continuity} to \eqref{eq:fx_constit1} results in
\begin{align}
&\sum_{i\in\partial_+j}\int_{\overline{\Omega}_{ij}}\partial_x\phi_{ij}(t,x_{ij}) +\sum_{k\in\partial_-j}\int_{\underline{\Omega}_{ij}}\partial_x\phi_{jk}(t,x_{jk}) \label{eq:pt_constit1}  \\
&\,\, =  -\!\!\sum_{i\in\partial_+j}\int_{\overline{\Omega}_{ij}}\!\!\partial_t \rho_{ij}(t,x_{ij}) - \!\!\sum_{k\in\partial_-j}\int_{\underline{\Omega}_{ij}}\!\!\partial_t \rho_{jk}(t,x_{jk})   \label{eq:pt_constit2}\\
&\,\,  \approx -\sum_{i\in\partial_+j}\frac{\eP}{2}\dot{\overline{\rho}}_{ij} -\sum_{k\in\partial_-j}\frac{\eP}{2}\dot{\underline{\rho}}_{jk}
\label{eq:pt_constit3} \\
&\,\, =  -\sum_{i\in\partial_+j}\frac{\eP}{2}\ddx{t}(\overline{\alpha}_{ij}\rho_{j}) -\sum_{k\in\partial_-j}\frac{\eP}{2}\ddx{t}(\underline{\alpha}_{jk}\rho_{j})
\label{eq:pt_constit4} \\
&\,\, =  -\frac{\eP}{2}\!\!\sum_{i\in\partial_+j}\!(\dot{\overline{\alpha}}_{ij}\rho_{j} + \overline{\alpha}_{ij}\dot{\rho}_{j}) -\frac{\eP}{2}\!\!\sum_{k\in\partial_-j}\!(\dot{\underline{\alpha}}_{jk}\rho_{j}+ \underline{\alpha}_{jk}\dot{\rho}_{j})
\nonumber \\
&\,\, = -\frac{\eP}{2}(\dot{\alpha}_j\rho_j+ \alpha_j \dot{\rho}_j), \label{eq:pt_constit6}
\end{align}
where $\alpha_j$ denotes aggregated actuation at node $j\in\hat{\cV}_{\eP}$,
\begin{align}
\alpha_j(t)= \sum_{i\in\partial_+j} \overline{\alpha}_{ij}(t) + \sum_{k\in\partial_-j} \underline{\alpha}_{jk}(t).
\label{eq:alpha_def}
\end{align}
The approximation in \eqref{eq:pt_constit3} is made by assuming sufficient network refinement \eqref{eq:pres_rel}, and the nodal density relations \eqref{eq:in_pressure_comp} are substituted into \eqref{eq:pt_constit3} to obtain \eqref{eq:pt_constit4}.  We have established equality of \eqref{eq:fx_constit3} and \eqref{eq:pt_constit6}, so solving for $\rho_j$ yields the discretized nodal mass conservation dynamics
\begin{align}
\!\!\!\!\!\! \dot{\rho}_j \!=\! \dS \frac{2}{\eP\alpha_j}\bq{\sum_{i\in\partial_+j}\!\!\vphi_{ij} \!\!  - \!\!\!\dS\sum_{k\in\partial_-j}\!\!\vphi_{jk} \!+ q_j} \! - \! \frac{\dot{\alpha}_j}{\alpha_j}\rho_j, \, \fA j\in\hat{\cV}_{\eP}.  \label{eq:disc_mass_balance}
\end{align}
The units on both sides of \eqref{eq:disc_mass_balance} are equal, e.g. $kg/m^3/s$.

Next, we approximate the dissipation equation \eqref{eq:in_dissipation_eq} by evaluating the spatial gradient with a finite difference
\begin{align}
\partial_{x}\rho_{ij}(t,x_{ij}) \approx \frac{1}{\eP}(\overline{\rho}_{ij}-\underline{\rho}_{ij}) = \frac{1}{\eP}(\overline{\alpha}_{ij}\rho_j-\underline{\alpha}_{ij}\rho_i), \label{eq:px_grad}
\end{align}
accounting for endpoint actuators as shown at right in Figure \ref{fig:netconstit}.  Applying \eqref{eq:px_grad} to approximate \eqref{eq:in_dissipation_eq} at nodes yields
\begin{align}
\!\!\!\!\! \vphi_{ij} \!\!&=\!\! \dS -f_{\mu(ij)}\!\bp{\! t,\overline{\alpha}_{ij}\rho_{j}, \!\frac{1}{\eP}(\overline{\alpha}_{ij}\rho_j \! -\underline{\alpha}_{ij}\rho_i) \!}, \,\,\,\,\,\, \fA i\in\partial_+j,  \label{eq:disc_diss_eq1} \\
\!\!\!\!\! \vphi_{jk} \!\!&=\!\! \dS -f_{\mu(jk)}\!\bp{\! t,\underline{\alpha}_{jk}\rho_{j}, \!\frac{1}{\eP}(\overline{\alpha}_{jk}\rho_k \! -\underline{\alpha}_{jk}\rho_j) \!}, \, \fA  k\in\partial_-j,  \label{eq:disc_diss_eq2}
\end{align}
where $(i,j)$ and $(j,k)$ are used for incoming and outgoing edges at node $j$, respectively.  The density $\overline{\alpha}_{ij}\rho_j$ (resp. $\underline{\alpha}_{jk}\rho_j$) is used in the second argument of $f_{\mu(ij)}$ (resp. $f_{\mu(jk)}$), rather than an average of that and $\underline{\alpha}_{ij}\rho_i$ (resp. $\overline{\alpha}_{jk}\rho_k$), because of the assumption in \eqref{eq:pres_rel}.  Substituting \eqref{eq:disc_diss_eq1}-\eqref{eq:disc_diss_eq2} into \eqref{eq:disc_mass_balance} produces the purely nodal dynamics
\begin{align}
 \dS \dot{\rho}_j & = \dS \frac{2}{\eP\alpha_j}\sum_{k\in\partial_-j} f_{\mu(jk)}\bp{t,\underline{\alpha}_{jk}\rho_j, \frac{1}{\eP}(\overline{\alpha}_{jk}\rho_k-\underline{\alpha}_{jk}\rho_j)} \nonumber \\ & \qquad -  \frac{2}{\eP\alpha_j}\dS\sum_{i\in\partial_+j} f_{\mu(ij)}\bp{t,\overline{\alpha}_{ij}\rho_j, \frac{1}{\eP}(\overline{\alpha}_{ij}\rho_j-\underline{\alpha}_{ij}\rho_i)}\nonumber \\ & \qquad + \frac{2}{\eP\alpha_j}q_j - \frac{\dot{\alpha}_j}{\alpha_j}\rho_j , \quad \fA j\in\hat{\cV}_{\eP}.  \label{disceq3}
\end{align}

\begin{rem}{Regularity Assumptions.} \label{rem:mol} First, we note that the ODE system \eqref{disceq3} is defined on the nodes $\hat{\cV}_\eP$ of the $\eP$-refined graph $\cG_\eP$.  We assume that this discretization scheme for the PDE system defined by \eqref{eq:in_continuity}-\eqref{eq:in_dissipation_eq} with \eqref{eq:in_nodal_continuity}-\eqref{eq:in_pressure_comp} is convergent and stable in the sense of a method of lines (MOL) solution.  That is, the distance between solutions to \eqref{disceq3} and the classical solution to the PDE system at locations corresponding to refined network nodes will converge point-wise to as $\eP\to 0$.
\end{rem}

We derive nodal density dynamics to enable investigation of monotonicity properties of the network flow PDE system with respect to a finite collection of parameter functions $q_i$.  We propose that if the assumptions of Remarks \ref{rem:wellposed} and \ref{rem:mol} hold, then such monotonicity properties derived for the ODE system will hold for the PDE system.  In this paper, we do not attempt to prove well-posedness, regularity, and convergence of the approximation. Instead our intention is to derive a property assuming that such conditions hold, and which can be used for simplification of otherwise intractable optimal control problems subject to parameter uncertainty.  More specifically, we demonstrate a monotonicity property of the spatially discretized nodal ODE dynamics \eqref{disceq3}, and which is sufficient for efficient formulation of the associated robust optimal control problem, as described in Section \ref{sec:application}.

\section{Monotonicity of Nodal Dynamics} \label{sec:monotonicity}

The monotonicity property that we will derive for the system \eqref{disceq3} states that nodal density can only increase with increasing injection.  Therefore, we consider monotonicity with respect to a time-varying parameter function rather than the control input, in contrast to the monotone control systems literature \cite{angeli03,como13a}.  This requires several definitions in a control system setting.

\begin{df}[Parameterized control system] Consider
\begin{align} \label{sys0}
\dot{x}=g(x,u,p), \qquad x(0)=y
\end{align}
with state $x(t)\in\cX\subset\bR^n$, control vector $u(t)\in\cU\subset\bR^m$, and parameter vector $p(t)\in\cP\subset\bR^p$ where $g$ is Lipschitz continuous and $\cX$, $\cU$, and $\cP$ are closed and convex.
\end{df}

\begin{df}[Monotone parameterized control system] \label{monotonesysdef} The control system \eqref{sys0} is \emph{monotone parameterized} with respect to $p(t)$ if, for all $t\geq 0$ $y_1,y_2\in\cX$, $u(t):(0,\infty)\to\cU$, and piecewise-continuous functions $p_1(t), p_2(t):(0,\infty)\to\cP$, the orderings $y_1\leq y_2$ and $p_1(s)\leq p_2(s)$ $\fA s\in[0,t]$ imply that $x_1(t)\leq x_2(t)$.  Here the inequalities for vectors are meant componentwise (i.e., $y\leq z$ means that $y_i\leq z_i$ for all $i=1,\ldots,n$), and $x_j(t)$, for $j=1,2$, stands for the solution to \eqref{sys0} with initial condition $y_j$, control input $u(t)$, and parameter vector $p_j(t)$.
\end{df}

\begin{df}{Non-negative matrix.} A square matrix $B\in\bR^{n\times n}$ is called \emph{non-negative} if all its entries are non-negative.
\end{df}

\begin{df}{Metzler matrix.} A square matrix $A\in\bR^{n\times n}$ is called \emph{Metzler} if all its off-diagonal entries are non-negative, i.e., $A_{ij}\geq 0$ for all $i\neq j\in\{1,\ldots, n\}$.
\end{df}

The main result in this section depends on the following well-known monotone dynamical systems theorem.

\begin{thm}[Dynamic Monotonicity Conditions] \label{monotonesystheo} The dynamical system \eqref{sys0} is a monotone parameterized control system if and only if $\nabla_x g$ is Metzler and $\nabla_p g$ is non-negative almost everywhere in $\bR^n\times\bR^m$ for all $u(t):(0,\infty)\to\cU$.
\end{thm}

\pf A standard application of the Kamke-M\"uller conditions \cite{kamke32,angeli03,hirsch05}. \hfill \boX



\begin{prop}[Monotonicity of Nodal Flow Dynamics] \label{propgasmonotone} The nodal network flow dynamics in \eqref{disceq3} are monotone parameterized with respect to nodal commodity injections $q=(q_1,\ldots,q_V)$ given positive compression ratio functions $\underline{\alpha}_{ij},\overline{\alpha}_{ij}\in C_+^k([0,T])$ for $k\geq 1$ if
the dissipation function $f_{ij}(t,u,v)$ is differentiable and increasing in its last argument for all original network edges, i.e.,
\begin{align}
\ppx{v} f_{ij}(t,u,v)>0, \quad (i,j)\in\cE. \label{eq:in_increasing_dissipation}
\end{align}
\end{prop}
\vspace{1ex}

\pf The derivatives of $\dot{\rho}_j$ with respect to $\rho_i$ (for $i\in\partial_+j$), $\rho_k$ (for $k\in\partial_-j$), and $q_m$ (for any $m\in\{1,\ldots,V\}$), respectively, are
\begin{align}
 \!\!\dS \pxpy{\dot{\rho}_j}{\rho_i} &= \dS \frac{2\overline{\alpha}_{ij}}{\eP^2\alpha_j}h_{\mu(ij)}\bp{ t,\overline{\alpha}_{ij}\rho_j, \frac{1}{\eP}(\overline{\alpha}_{ij}\rho_j-\underline{\alpha}_{ij}\rho_i)} \label{eq:monderiv1} \\
 \!\!\dS \pxpy{\dot{\rho}_j}{\rho_k} &= \dS \frac{2\underline{\alpha}_{jk}}{\eP^2\alpha_j}\dS h_{\mu(ij)}\bp{ t,\underline{\alpha}_{jk}\rho_j, \frac{1}{\eP}(\overline{\alpha}_{jk}\rho_k-\underline{\alpha}_{jk}\rho_j)} \label{eq:monderiv2} \\
 \!\!  \dS \pxpy{\dot{\rho}_j}{q_m}  &= \dS \frac{2}{\eP\alpha_j}\chi_{jm} \label{eq:monderiv3}
\end{align}
where $\chi_{jm}=1$ if $j=m$ and $\chi_{jm}=0$ if $j\neq m$.
Derivatives of the dissipation functions satisfy $h_{\mu(ij)}(t,u,v)\triangleq\ppx{v}f_{\mu(ij)}(t,u,v)>0$ for all $(i,j)\in\hat{\cE}_{\eP}$, which is inherited from $\ppx{v}f_{ij}(t,u,v)>0$ for $(i,j)\in\cE$, as stipulated in \eqref{eq:in_increasing_dissipation}.  In addition, $\overline{\alpha}_{ij}(t),\underline{\alpha}_{ij}(t)>0$ for all $(i,j)\in\hat{\cE}_{\eP}$ and $\alpha_j(t)>0$ for all $j\in\hat{\cV}_{\eP}$ at all $t>0$ as well. Therefore \eqref{eq:monderiv1}-\eqref{eq:monderiv3} are strictly non-negative, so the gradient of the dynamics with respect to the state $\rho=(\rho_1,\ldots,\rho_V)$ is Metzler, and the gradient with respect to $q=(q_1,\ldots,q_V)$ is non-negative.    Thus the conditions of Theorem \ref{monotonesystheo} hold. \hfill \boX

\begin{rem} \label{remgasmonotone} Consider two initial nodal density vectors $\rho^1(0),\rho^2(0)\in\bR_+^V$, two vectors of continuous nodal injection functions $q^1,q^2:[0,t]\to\bR^V$, and a collection of positive compression ratio functions $\underline{\alpha}_{ij},\overline{\alpha}_{ij}:[0,T]\to\bR_+$ for $(i,j)\in\cC$.  Then if $\rho(0)^1\leq \rho(0)^2$ and $q^1(s)\leq q^2(s)$ for all $s\in[0,t]$, it follows that $\rho^1(t)\leq \rho^2(t)$, where inequalities for vectors are meant componentwise.  Here $\rho^k(t)$, for $k=1,2$, denote solutions to \eqref{disceq3} for all $j\in\cV$ and with initial condition $\rho^k(0)$ and injection vector $q^k(t)$.  This furthermore implies that for initial densities $\rho(0)\in\bR_+^V$ and nodal injection functions $q:[0,t]\to\bR^V$ that satisfy $\rho^1(0)\leq \rho(0) \leq \rho^2(0)$ and $q^1(s)\leq q(s) \leq q^2(s)$ $\fA \, s\in[0,t]$, then the solution to \eqref{disceq3} will satisfy $\rho^1(t)\leq \rho(t) \leq \rho^2(t)$ componentwise.  This property holds when actuator functions $\underline{\alpha}_{ij},\overline{\alpha}_{ij}$ are provided as state-independent parameters.
\end{rem}

\section{Application to Robust Optimal Control} \label{sec:application}

The results of Sections \ref{sec:discretization} and \ref{sec:monotonicity} enable us to establish important properties and robust optimal control formulations for dynamic dissipative flows on networks.

\subsection{Local Feedback Control can Preserve Monotonicity }

Proposition \ref{propgasmonotone} can be extended in a straightforward manner to produce sufficient conditions on local feedback policies for commodity flow actuation to preserve monotonicity.

\begin{cor}[Monotonicity-Preserving Local Feedback] \label{propfeedback} The nodal dynamics in \eqref{disceq3} are monotone parameterized with respect to nodal injections $q=(q_1,\ldots,q_V)$ given local feedback policies $\overline{\alpha}_{ij}(t)=\overline{k}_{ij}(\rho_j(t))$ for $i\in\partial_+j$ and $\underline{\alpha}_{jk}(t)=\underline{k}_{jk}(\rho_j(t))$ for $k\in\partial_-j$ if the functions $\underline{k}_{ij}(v)$ and $\overline{k}_{ij}(v)$ satisfy
$v\ddx{v}\underline{k}_{ij}(v)+\underline{k}_{ij}(v)>0$ and $v\ddx{v}\overline{k}_{ij}(v)+\overline{k}_{ij}(v)>0$ for all $(i,j)\in \cC$ and $v\in\bR_+$.
\end{cor}

\pf Define the function
\begin{align}
k_j(v)= \sum_{i\in\partial_+j} \overline{k}_{ij}(v) + \sum_{k\in\partial_-j} \underline{k}_{jk}(v),
\label{eq:kay_def}
\end{align}
with $k'(v)\equiv\ddx{v}k(v)$, so that under the local feedback policies the final term in the right hand side of the nodal dynamics \eqref{disceq3} becomes $k'_j(\rho_j)\dot{\rho}_j\rho_j/k_j(\rho_j)$.
Thus, solving for $\dot{\rho}_j$ yields the nodal dynamics
\begin{align}
 \dS \dot{\rho}_j & = \dS \frac{2}{\eP r_j}\sum_{k\in\partial_-j} f_{\mu(jk)}\bp{t,\underline{k}_{jk}\rho_j, \frac{1}{\eP}(\overline{k}_{jk}\rho_k-\underline{k}_{jk}\rho_j)} \nonumber \\ & \qquad - \frac{2}{\eP r_j}\dS\sum_{i\in\partial_+j} f_{\mu(ij)}\bp{t,\overline{k}_{ij}\rho_j, \frac{1}{\eP}(\overline{k}_{ij}\rho_j-\underline{k}_{ij}\rho_i)} \nonumber \\ & \qquad + \frac{2}{\eP r_j}q_j, \quad \fA j\in\hat{\cV}_{\eP},  \label{eq:disc_feedback1}
\end{align}
where $r_j(\rho_j)=k_j(\rho_j)+k_j'(\rho_j)\rho_j$. Following the procedure in Proposition \ref{propgasmonotone}, the gradient of \eqref{eq:disc_feedback1} with respect to $q=(q_1,\ldots,q_V)$ is non-negative and the gradient with respect to  $\rho=(\rho_1,\ldots,\rho_V)$ is Metzler when $r_j(\rho_j)>0$ $\fA$ $j\in\hat{\cV}_{\eP}$.  This holds given the above assumptions on $\underline{k}_{ij}$ and $\overline{k}_{ij}$ . \boX

\begin{figure}[t] 
 \centering{
\includegraphics[width=.45\linewidth]{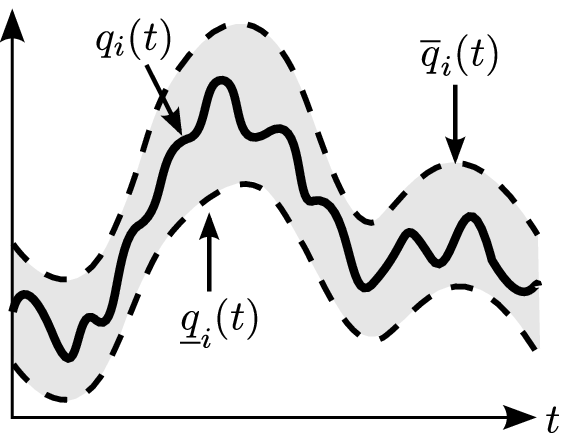} \quad \includegraphics[width=.45\linewidth]{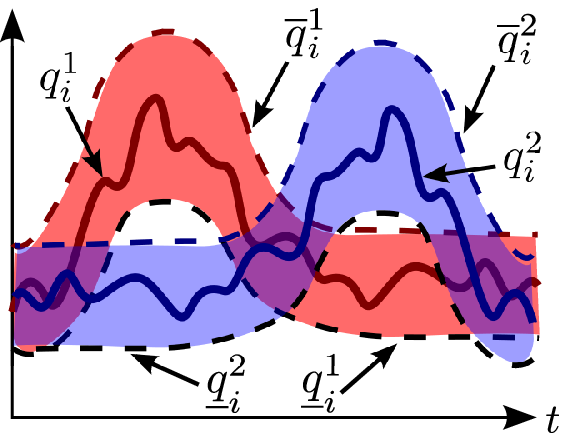}  \caption{Example single (left) or multiple (right) injection uncertainty regions.} \label{fig:monoton1}}
\end{figure}

\subsection{Robust Optimal Control Formulation}

Proposition \ref{propgasmonotone} also enables compact formulations that greatly simplify optimal control problems for dynamic commodity flow networks in which global information is available in advance.  We consider the case when the commodity injections $q_i(t)$ are known for an entire time interval $[0,T]$ so that the control protocols $\underline{\alpha}_{ij}$ and $\overline{\alpha}_{ij}$ can be determined predictively using the dynamic model \eqref{eq:in_continuity}-\eqref{eq:in_dissipation_eq} with \eqref{eq:in_nodal_continuity}-\eqref{eq:in_pressure_comp}.  This type of optimal control problem appears in dynamic commodity transport applications, and the monotonicity property enables tractable extension to situations with uncertain injections $q_i$. 

Observe that the discretized network flow system \eqref{disceq3} can be written in the form of \eqref{sys0} as $\dot{\rho}(t) = g(\rho(t),\alpha(t),q(t))$, where $\rho(t)\in\bR_+^V$ is the state, $\alpha(t)\in\bR_+^C$ is a control vector where $C=|\cC|$ is the number of compression actuators, and $q(t)\in\bR^V$ is the parameter vector of nodal commodity injections. Consider the optimal control problem given by
\begin{align}
	\min_{\alpha}\ \ & J(\rho,\alpha)=\int_0^T \mathcal{L}(t,\rho(t),\alpha(t))dt, \label{eq:ocp0a} \\
	{\rm s.t.}\ \ & \dot{\rho}(t) = g(\rho(t),\alpha(t),q(t)), \label{eq:ocp0b} \\
	& \rho_{\min} \leq \rho(t)\leq \rho_{\max}. \label{eq:ocp0c}
\end{align}
Here $\mathcal{L}\in C^k$ is in the space $C^k$ of continuous functions with $k\in\bN$ classical derivatives, and the dynamic constraints $g\in C_V^{k-1}$ are in the space $C_V^{k-1}$ of $V$-vector valued $C^{k-1}$ functions, with respect to the state, $\rho(t)\in\bR^V$, and control, $\alpha(t)\in\bR^C$.  The admissible set for controls $\alpha$ includes the piecewise $C_m^k$ functions on $[0,T]$.  The values $\rho_{\min}$ and $\rho_{\max}$ are minimum and maximum commodity densities at network nodes, and the vector inequalities \eqref{eq:ocp0c} are entry-wise.

Suppose that each commodity injection $q_i\in C^{k}([0,T])$ is an uncertain function $q_i(t)\in[\underline{q}_i(t),\overline{q}_i(t)]$ (or even $q_i(t)\in[\underline{q}_i^1(t),\overline{q}_i^1(t)]\cup[\underline{q}_i^2(t),\overline{q}_i^2(t)]$) for all $t\in[0,T]$ as shown in Figure \ref{fig:monoton1}.  Then the optimal control problem \eqref{eq:ocp0a}-\eqref{eq:ocp0c} becomes very challenging because the dynamic constraints \eqref{eq:ocp0b} are repeated for an uncountable continuum of possible functions $q_i\in C^{k}([0,T])$, and any sampling approach quickly becomes intractable.  Proposition \ref{propgasmonotone} enables the reformulation of \eqref{eq:ocp0a}-\eqref{eq:ocp0c} into the expanded problem
\begin{align}
	\min_{\alpha(t)}\ \ & J(\rho,\alpha)=\int_0^T ({\mathcal{L}}(t,\widetilde{\rho}(t),\alpha(t)) dt, \label{eq:ocp1a} \\
	{\rm s.t.}\ \ & \dot{\rho}(t) = g(\rho(t),\alpha(t),q(t)), \label{eq:ocp1b} \\
	&\dot{\underline{\rho}}(t) = g(\underline{\rho}(t),\alpha(t),\underline{q}(t)), \label{eq:ocp1b1} \\
& \dot{\overline{\rho}}(t) = g(\overline{\rho}(t),\alpha(t),\overline{q}(t)), \label{eq:ocp1b2} \\
	& \rho_{\min} \leq \underline{\rho}(t)\leq \rho_{\max}. \label{eq:ocp1c1} \\
& \rho_{\min} \leq \overline{\rho}(t)\leq \rho_{\max}. \label{eq:ocp1c2}
\end{align}
In the above formulation, $\widetilde{\rho}(t)$ can be one of $\underline{\rho}(t), \overline{\rho}(t)$ and $\rho(t)$ depending on the objective function we choose to optimize. Choosing
$\widetilde{\rho}(t) = \rho(t)$ corresponds to minimizing the nominal operational cost associated with the nominal injection values $q(t)$. Furthermore, if we assume that the cost functional
$\mathcal{L}(t,\rho(t),\alpha(t))$ is monotone with respect to $\rho(t)$, then we can optimize a worst case min-max objective by substituting $\widetilde{\rho}(t)$ to be one of
$\overline{\rho}(t)$ or $\underline{\rho}(t)$ depending on the sign of the monotonicity of $\mathcal{L}$ with respect to its second entry. The min-max objective is defined as
\begin{align}
	J_{mm}(\rho,\alpha) =  \max_{q(t) \in [\underline{q}(t), \overline{q}(t)]} & \int_0^T {\mathcal{L}}(t,{\rho}(t),\alpha(t)) dt, \label{eq:ocp2a} \\
	{\rm s.t.}\ \ & \dot{\rho}(t) = g(\rho(t),\alpha(t),q(t)). \label{eq:ocp2b}
\end{align}
Because by Remark \ref{remgasmonotone} $\rho(t)$ increases monotonically with respect to $q(t)$, and by assumption $\mathcal{L}$ is monotone with respect to $\rho(t)$, the maximum in \eqref{eq:ocp2a} is obtained by substituting
$q(t)$ to be one of $\underline{q}(t)$ or $\overline{q}(t)$, which justifies the choice of $\widetilde{\rho}(t)$ described above.
Observe that we have enforced upper and lower feasibility bounds in \eqref{eq:ocp1c1} and \eqref{eq:ocp1c2} only on $\overline{\rho}(t)$ and $\underline{\rho}(t)$.
Suppose that a control protocol $\alpha(t)$ is found to minimize the running cost \eqref{eq:ocp1a} to satisfy the dynamic and box constraints \eqref{eq:ocp1b1} and \eqref{eq:ocp1c1} for the maximum injection vector $\overline{q}(t)$ and the same constraints \eqref{eq:ocp1b2} and \eqref{eq:ocp1c2} for the minimum injection $\underline{q}(t)$.  As described in Remark \ref{remgasmonotone}, this is sufficient for $\alpha(t)$ to also be feasible for any injection profile that satisfies $\underline{q}(t)\leq q(t) \leq \overline{q}(t)$.  This key property is of significance to creating tractable optimal control formulations for uncertain dynamic dissipative network flows.

\section{Conclusions} \label{sec:conc}

We have derived monotonicity properties for a class of actuated dynamic flows on networks described by dissipative partial differential equation systems with nodal consistency conditions, which are approximated by discretization in space.  The result shows that discretization can preserve the desirable monotonicity property of such network flows, which was recently established \emph{ab initio} \cite{misra16mtns}.  Specifically, this states that under certain conditions, commodity density anywhere in the network can only increase monotonically when any commodity injection is increased.   We also characterized the conditions on local feedback policies for actuators throughout the network for which monotonicity is maintained. In addition, the results enable compact formulations that greatly simplify canonical optimal control problems involving uncertain dynamic commodity flows over networks, in which the solutions must be robust with respect to commodity withdrawal uncertainty within maximum and minimum time-dependent limits.  
These results are expected to find applications in large-scale systems for transportation of energy \cite{budivsic12,zlotnik15cdc,vuffray15cdc,zlotnik15dscc}, vehicles \cite{dapice08,ma04,tjandra03phd,gottlich05,gugat05,herty03,como10,como15}, and information \cite{como13a,como13b,intanagonwiwat02}.

\section*{Acknowledgements}
We thank Giacomo Como for valuable discussions.  This work was carried out under the auspices of the National Nuclear Security Administration of the U.S. Department of Energy at Los Alamos National Laboratory under Contract No. DE-AC52-06NA25396, and was partially supported by DTRA Basic Research Project \#10027-13399 and by the Advanced Grid Modeling Program in the U.S. Department of Energy Office of Electricity.





\bibliographystyle{unsrt}
\bibliography{gas_master,monotonicity_master}

\end{document}